# A Distributionally Robust Self-Scheduling Under Price Uncertainty Based on CVaR

Linfeng Yang, *Member, IEEE*, Ying Yang, Guo Chen, *Member, IEEE*, Zhaoyang Dong, *Fellow, IEEE*

*Abstract*—To ensure a successful bid while maximizing of profits, generation companies (GENCOs) need a self-scheduling strategy that can cope with a variety of scenarios. So distributionally robust optimization (DRO) is a good choice because that it can provide an adjustable self-scheduling strategy for GENCOs in the uncertain environment, which can well balance robustness and economics compared to strategies derived from robust optimization (RO) and stochastic programming (SO). In this paper, a novel moment-based DRO model with conditional value-at-risk (CVaR) is proposed to solve the self-scheduling problem under electricity price uncertainty. Such DRO models are usually translated into semi-definite programming (SDP) for solution, however, solving large-scale SDP needs a lot of computational time and resources. For this shortcoming, two effective approximate models are proposed: one approximate model based on vector splitting and another based on alternate direction multiplier method (ADMM), both can greatly reduce the calculation time and resources, and the second approximate model only needs the information of the current area in each step of the solution and thus information private is guaranteed. Simulations of three IEEE test systems are conducted to demonstrate the correctness and effectiveness of the proposed DRO model and two approximate models.

*Index Terms*—self-scheduling, distributionally robust optimization, price uncertainty, CVaR, ADMM.

## NOMENCLATURE

Operator:
$[\cdot]^+$    $\max(0,\cdot)$.
$|\cdot|$    Number of elements in "·", if "·" is a set.
$[N]$    Set $\{1, \dots, N\}$.
$\Re^{\cdot}$    $\Re$ is real number set, "·" is dimension of set.

Indices:
$i, j$    Index for unit/bus.
$t$    Index for time periods.

Constants:
$S^\text{B}$    Set of buses, $|S^\text{B}| = N$ and $S^\text{B} = [N]$.
$S^\text{G}$    Set of buses with units. $|S^\text{G}| = G$ and $S^\text{G} \subseteq S^\text{B}$.
$S^\text{TP}$    Set of time periods, and $|S^\text{TP}| = T$.
$S^\text{L}$    Set of lines, and $|S^\text{L}| = L$.
$S^\text{ref}$    Set of reference bus.
$a_i, b_i, c_i$    Coefficients of the quadratic production cost function of unit $i$.
$F_{ij}^\text{max}$    Maximum power flow limit on line $ij$.
$X_{ij}$    Reactance of line $ij$.
$\underline{P}_i^\text{G}, \bar{P}_i^\text{G}$    Minimum and maximum power output of unit $i$.
$P_{i,t}^\text{D}$    System load demand on bus $i$ in time periods $t$.
$P_i^\text{G,up}, P_i^\text{G,down}$    Ramp up and down limit of unit $i$.
$\boldsymbol{B}$    Network admittance matrix ($\Re^{N \times N}$).
$B_{ij}$    Element of $\boldsymbol{B}$ in the $i^\text{th}$ row and $j^\text{th}$ col.

Variables:
$P_{i,t}^\text{G}$    Power output of unit on bus $i$ in time periods $t$.
$\theta_{i,t}$    Phase angle on bus $i$ in time periods $t$.
$\lambda_{i,t}$    The electricity price of unit $i$ (depending on the bus where the unit is connected) in period $t$.
$\boldsymbol{P}, \boldsymbol{\theta}, \boldsymbol{\lambda}$    Represent $\left[P_{[i],[t]}^\text{G}\right]^\text{T}$, $\left[\theta_{[i],[t]}\right]^\text{T}$, and $\left[\lambda_{[i],[t]}\right]^\text{T}$, respectively.

## I. INTRODUCTION

IN the day-ahead or real-time power energy market, for the sake of maximize the generation profit, generation companies (GENCOs) usually need to generate bidding schemes according to the self-scheduling strategy, which is generally obtained by solving the profit-maximizing optimal power flow model based on the locational marginal prices(LMPs) [1], [2]. For the bidding scheme to be accepted by independents system operators (ISOs), an appropriate self-scheduling strategy should to be adopted. Conejo et al. [2] provides a framework for solving the optimal bidding scheme of price-taker producer based on price forecasting. And a deterministic solution for risk measurement and management of GENCOs is presented in [3]. Further, Yamin and Shahidehpour [4] considered the self-scheduling problems in formulations that account for profit and risk simultaneously. However, in this kind of real-time self-scheduling problem, it is difficult to predict the price accurately, which leads to uncertainty. So, a methodology has been proposed for risk management based on fuzzy numbers to model price uncertainty in [5], [6]. At the same time, with the development of power energy market, more and more uncertainty variables are introduced to make the model better describe the actual scenario. In this uncertain environment, how to obtain an appropriate self-scheduling strategy has attracted the attention of scholars and engineer [7]. In the past decades, there are two main methods to solve this kind of problems with uncertainty variables: robust optimization (RO) and stochastic programming (SO) [8]. RO constructs an uncertainty set in which all uncertainty variables fall, by optimizing the worst case in the uncertainty set, the robustness of the scheduling scheme is guaranteed [9], [10]. While SO uses methods such as average sampling approximation and random approximation to obtained a distribution of uncertainty variable and always as-

The work of L.F. Yang and Y. Yang were supported by the Natural Science Foundation of China (51767003), and Natural Science Foundation of Guangxi (2020GXNSFAA297173, 2020GXNSFDA238017). The work of G Chen and ZY Dong are partially supported by UNSW Digital Grid Futures Institute seed fund.

L.F. Yang and Y. Yang are with the School of Computer Electronics and Information, Guangxi University, Nanning 530004, China, and L.F. Yang is also with the Guangxi Key Laboratory of Multimedia Communication and Network Technology, Guangxi University. (e-mail: ylf@gxu.edu.cn; 907803678@qq.com).

G. Chen (the corresponding author) and Z.Y. Dong are with the School of Electrical Engineering and Tele-communications, The University of NSW, Sydney, NSW 2052, Australia. (e-mail: andyguochen@gmail.com; zydong@ieee.org)

sumes that the distribution is real [11]-[13]. However, the above two methods have some drawbacks, RO provides security in worst case, but it always generates over-conservative self-scheduling strategy in practice because it ignores the potential distribution information of uncertainty variable, and SO always assumes that the distribution of uncertainty variable is known, but the assumed distribution in practical application may not accurate, and resulting in the generation of wrong strategy [8], [14].

For this reason, a method called distributionally robust optimization (DRO) has gradually attracted considerable attention in the academic field. The DRO was first proposed by Zackova in [15], and it can be traced back to Scarf's [16] optimal single product newsboy problem under the unknown distribution of known mean and variance. DRO is proposed to balance the conservatism of RO and the particularity of SO, its optimization decision is to find the worst probability distribution in the possible distribution family and then optimize it, so that all probability distributions are immunized. DRO lies all possible probability distributions in a set called ambiguity set, which is data-driven and assumed to contain the true distribution, instead of constructing a deterministic uncertainty set [17]. According to the method of constructing ambiguity set, DRO can be divided into moment-based DRO [18]-[23] and distance-based [8], [17] DRO. After years of research and development, DRO, especially the moment-based DRO, is used in many fields such as energy system and unit commitment, etc. An auxiliary random variable is introduced to describe the distribution characteristics of wind power generation in [18], then a new two-stage DRO model for UC scheduling under wind uncertainty has been presented. [21] based on first- and second- moment information, a real-time DRO power dispatch model is established to obtain a scheme, which includes economic dispatch and corrective control. [22] presents a distributionally robust chance constrained planning method for the integrated heat and electricity system. However, such DRO models are not easy to solve, the study in [24] proposed an internal approximation method to solve the blockage management model of DRO, its solution strategy is easy to handle but may be not accurate enough. Then [19] proposed that the DRO model can be reconstructed into semi-definite programming (SDP) or second-order cone programming (SOCP) based on some mild assumptions, and it is further shown to be solvable in time polynomial. Overall, DRO uses the method of constructing ambiguity set to hedge against the inaccuracy of the probability distribution of uncertainty variable. Different ambiguity set construction methods produce different DRO model with different balance between robustness and economy [14]. The above characteristics make the DRO consider the potential distribution information of uncertainty variable without assuming the exact distribution of it, which well balance the advantages and disadvantages of RO and SO [25].

In above works, most of them are optimizing generation costs, and the potential risks caused by uncertainty are not well managed, therefore, in reference [26], in order to better characterize risks and carry out risk management, the value-at-risk (VaR) tools commonly used in the field of economics portfolio selection are usually introduced. In the self-scheduling problem, it specifically refers to the probability that the possible loss of the generated scheduling scheme under a normal fluctuation of uncertain parameters does not exceed a certain set value meets the confidence level. For avoid the deficiencies of VaR in excess loss control, the conditional value-at-risk (CVaR) should be introduced based on VaR [1], [20], [27].

This paper proposes a DRO model for self-scheduling problem, which take care of price uncertainty, decision making and risks management. For each GENCO, an adjustable self-scheduling strategy including total cost, potential risks, and units dispatch is provided, decision makers can adjust the weights of risks and benefits to obtain different strategy according to the actual situation. Compared with the fixed strategy derived from the RO and SO models, an adjustable strategy is more flexible to adapt to the changing conditions of reality. The major contributions of this paper are summarized as follows:

1) In the self-scheduling problem, to remedy the SO's specificity and RO's over-conservative while controlling risk, a novel moment-based DRO model with CVaR is developed for the GENCOs under electricity price uncertainty. It can well combine the advantages of RO and SO while avoiding their disadvantages, so it gets a self-scheduling strategy that considers both robustness and economy.

2) It is different from the traditional fixed self-scheduling strategy, an adjustable strategy is derived from our model, which provides a promising alternative for the decision makers. The decision maker adjusts the ratio of robustness and economy in the final strategy by adjusting the size of ambiguity set according to the actual scenario, the tighter the ambiguity, the more aggressive the strategy is.

3) The computational burden grows heavily with the size of the problem, so two approximate models are proposed, both of them can significantly reduce the computational time and resources required while maintaining accuracy, and the second approximate model adopts distributed algorithm and thus can further protect the privacy.

The rest of this paper is organized as follows: Section II presents the formulation of the self-scheduling problem. A DRO framework for uncertainty self-scheduling and two approximations of DR-CVaR model are proposed in Section III and Section IV, respectively. Simulations results for IEEE 6- , 30- and 118-bus test system of the proposed DRO model and two approximate models are given in Section V while the conclusions are drawn in Section VI.

II. MATHEMATICAL FORMULATION OF SELF-SCHEDULING

A. Objective Function

The objective function of the self-scheduling problem is $F(\boldsymbol{y}, \boldsymbol{\lambda})$, which is the income of selling power minus the total operation cost [1], [28]
$$F(\boldsymbol{y}, \boldsymbol{\lambda}) = \boldsymbol{P}^\mathrm{T}\boldsymbol{\lambda} - f(\boldsymbol{P}),$$
where $\boldsymbol{y} = [\boldsymbol{P}; \boldsymbol{\theta}]$ is decision vector which includes all the decision variables. The fuel cost of units ($f(\boldsymbol{P})$) has the form
$$f(\boldsymbol{P}) = \sum_{i \in S^\mathrm{G}} \sum_{t \in S^\mathrm{TP}} f_i(P_{i,t}^\mathrm{G}), \tag{1}$$

where the production cost for unit $i$ can be represented as a quadratic function $f_i(\cdot)$, which is
$$f_i(\cdot) = a_i + b_i(\cdot) + c_i(\cdot)^2. \tag{2}$$

### B. Constraints

$$0 \leq P_{i,t}^G - \sum_{j \in S_B} B_{ij}\theta_{j,t} \leq P_{i,t}^D, \forall t \in S^{TP}, \forall i \in S^B, \tag{3}$$

$$\left|\frac{\theta_{i,t} - \theta_{j,t}}{x_{ij}}\right| \leq F_{ij}^{\max}, \forall t \in S^{TP}, \forall i,j \in S^B, \forall ij \in S^L, \tag{4}$$

$$\theta_{i,t} = 0, \forall i \in S^{\text{ref}}, \forall t \in S^{TP}, \tag{5}$$

$$P_{i,t}^G = 0, \forall i \in S^B - S^G, \forall t \in S^{TP}, \tag{6}$$

$$\underline{P}_i^G \leq P_{i,t}^G \leq \overline{P}_i^G, \forall i \in S^G, \forall t \in S^{TP}, \tag{7}$$

$$P_{i,t}^G - P_{i,t-1}^G \leq P_i^{G,up}, \forall i \in S^G, \forall t \in S^{TP}, \tag{8}$$

$$P_{i,t-1}^G - P_{i,t}^G \leq P_i^{G,down}, \forall i \in S^G, \forall t \in S^{TP}, \tag{9}$$

where Eq. (3) represents the DC network model, this inequality constraint shows that the total power generated by GENCOs should be less than or equal to the forecasted system demand. It is to be emphasized that GENCOs is not responsible for supplying the system demand and which is the ISO's responsibility [28]. Eq. (4) and Eq. (7) are capacity limits of transmission lines and generators. Eq. (5) defines the reference bus. Eq. (6) represents the buses without generators. Eq. (8)-(9) enforce ramp up/down limits of individual generators.

Assume that $L$ is a given integer parameter, let $p_{i,l} = \underline{P}_i + l(\overline{P}_i - \underline{P}_i)/L$ and $l = 0,1,2,\dots,L-1$.

After replacing $f_i(P_{i,t}^G)$ in the objective function with a corresponding new variable $z_{i,t}$ and adding the following linear constraints to the formulation,

$$z_{i,t} \geq (2c_i p_{i,l} + b_i)P_{i,t} + a_i - c_i(p_{i,l})^2, \tag{10}$$

we obtain a linear approximation of self-scheduling problem. For the sake of discussion, we use $Y$ to denote the feasible set ((3)-(10)) of self-scheduling defined in this subsection and the objective of this problem is maximizing $F(\boldsymbol{y}, \boldsymbol{\lambda}) = \boldsymbol{P}^T \boldsymbol{\lambda} - \mathbf{1}^T \boldsymbol{z}$, where $\boldsymbol{z} = [z_{[i],[t]}]^T$.

## III. DISTRIBUTIONALLY ROBUST CVaR OPTIMIZATION FRAMEWORK FOR UNCERTAIN SELF-SCHEDULING

### A. CVaR Based Uncertain Self-Scheduling

CVaR [27] is a well-known risk measure that has been widely used in various energy management problems [1], [29]. In this paper, CVaR is employed to hedge the risk faced by the self-scheduling under electricity price uncertainty.

Let the loss being $L(\boldsymbol{y}, \boldsymbol{\lambda}) = -F(\boldsymbol{y}, \boldsymbol{\lambda})$. If the uncertain $\boldsymbol{\lambda}$ obeys a distribution $\mathbb{P}$, then the probability of generating loss not exceeding a threshold $\alpha$ is then given by

$$\Psi(\boldsymbol{y}, \alpha) = \int_{L(\boldsymbol{y},\boldsymbol{\lambda}) \leq \alpha} p(\boldsymbol{\lambda}) d\boldsymbol{\lambda}, \tag{11}$$

where $p(\boldsymbol{\lambda})$ is the probability density function to $\mathbb{P}$.

The $\beta$-VaR and $\beta$-CVaR values for the loss random variables associated with $\boldsymbol{y}$ and any specified probability level $\beta \in (0,1)$, will be denoted by $\alpha_\beta(\boldsymbol{y})$ and $\phi_\beta(\boldsymbol{y})$. They are given by

$$\alpha_\beta(\boldsymbol{y}) = \min\{\alpha \in \Re : \Psi(\boldsymbol{y}, \alpha) \geq \beta\}, \tag{12}$$

$$\phi_\beta(\boldsymbol{y}) = \int_{L(\boldsymbol{y},\boldsymbol{\lambda}) \geq \alpha_\beta(\boldsymbol{y})} L(\boldsymbol{y}, \boldsymbol{\lambda}) \frac{p(\boldsymbol{\lambda})}{1-\beta} d\boldsymbol{\lambda}. \tag{13}$$

Then, $\alpha_\beta(\boldsymbol{y})$ is the lowest amount $\alpha$ such that, with probability $\beta$, the loss will not exceed $\alpha$. Whereas $\phi_\beta(\boldsymbol{y})$ is the conditional expectation of losses above $\alpha_\beta(\boldsymbol{y})$.

When we minimizing $\beta$-CVaR, a low $\beta$-VaR can be obtained as well because that $\beta$-VaR is never more than $\beta$-CVaR. Therefore, CVaR will be used to describe risks in the model:

Let $h(\boldsymbol{y}, \alpha, \boldsymbol{\lambda}) = \max\left\{\alpha, \frac{1}{1-\beta}L(\boldsymbol{y}, \boldsymbol{\lambda}) - \frac{\beta\alpha}{1-\beta}\right\}$, we have

$$\mathbb{E}_{\mathbb{P}}(h(\boldsymbol{y}, \alpha, \boldsymbol{\lambda})) = \alpha + \frac{1}{1-\beta}\mathbb{E}_{\mathbb{P}}([L(\boldsymbol{y}, \boldsymbol{\lambda}) - \alpha]^+) =$$

$$\alpha + \frac{1}{1-\beta}\int_{\boldsymbol{\lambda} \in \Re^{NT}}[L(\boldsymbol{y}, \boldsymbol{\lambda}) - \alpha]^+ p(\boldsymbol{\lambda}) d\boldsymbol{\lambda}. \tag{14}$$

And then, the following relationships can be verified according to the Theorem 1 and Theorem 2 in [27]:

$$\alpha_\beta(\boldsymbol{y}) \in \underset{\alpha \in \Re}{\operatorname{argmin}} \mathbb{E}_{\mathbb{P}}(h(\boldsymbol{y}, \alpha, \boldsymbol{\lambda})), \tag{15}$$

$$\phi_\beta(\boldsymbol{y}) = \min_{\alpha \in \Re} \mathbb{E}_{\mathbb{P}}(h(\boldsymbol{y}, \alpha, \boldsymbol{\lambda})) = \mathbb{E}_{\mathbb{P}}(h(\boldsymbol{y}, \alpha_\beta(\boldsymbol{y}), \boldsymbol{\lambda})), \tag{16}$$

$$\min_{\boldsymbol{y} \in Y} \phi_\beta(\boldsymbol{y}) = \min_{(\boldsymbol{y},\alpha) \in Y \times \Re} \mathbb{E}_{\mathbb{P}}(h(\boldsymbol{y}, \alpha, \boldsymbol{\lambda})). \tag{17}$$

### B. Distributionally Robust CVaR Based Self-Scheduling

While CVaR is an interesting risk measure, it still requires the decision maker to commit to a distribution $\mathbb{P}$. This is a step that can be difficult to take in practice. To introduce CVaR into our model, using the results of distributionally robust optimizations (DRO) [30], a general form of distributionally robust version of CVaR can be derived [31], i.e., given that the distribution is known to lie in a distributional set $\mathcal{D}_{M1}$, our distributionally robust CVaR (DR-CVaR) model of self-scheduling is, (based on the right hand side of (17))

$$\min_{(\boldsymbol{y},\alpha) \in Y \times \Re} \left\{ \sup_{\mathbb{P} \in \mathcal{D}_{M1}} \{\mathbb{E}_{\mathbb{P}}(h(\boldsymbol{y}, \alpha, \boldsymbol{\lambda}))\} \right\}. \tag{18}$$

### C. An Ambiguity Set for Uncertain Electricity Price

The ambiguity set of distribution in DRO model provides a flexible framework to model uncertainty by allowing the modelers to incorporate partial information about the uncertainty, obtained from historical data knowledge [30].

In what follows, we will study the DR-CVaR self-scheduling model under a popular version of distributional set based on moment [31]:

$$\mathcal{D}_{M1}(\mathcal{S}, \boldsymbol{\mu}, \Sigma, \gamma_1, \gamma_2) =$$

$$\left\{ \mathbb{P} \in \mathcal{P}_0(\Re^{GT}) \middle| \begin{array}{l} \mathbb{P}(\boldsymbol{\lambda} \in \mathcal{S}) = 1 \\ (\mathbb{E}_{\mathbb{P}}[\boldsymbol{\lambda}] - \boldsymbol{\mu})^T \Sigma^{-1}(\mathbb{E}_{\mathbb{P}}[\boldsymbol{\lambda}] - \boldsymbol{\mu}) \leq \gamma_1 \\ \mathbb{E}_{\mathbb{P}}[(\boldsymbol{\lambda} - \boldsymbol{\mu})(\boldsymbol{\lambda} - \boldsymbol{\mu})^T] \preccurlyeq \gamma_2 \Sigma \end{array} \right\}, \tag{19}$$

where $\mathcal{P}_0(\Re^{GT})$ denotes the set of all probability distributions on $\Re^{GT}$. $\mathcal{S} \subseteq \Re^{GT}$ is any closed convex set known to contain the support of $\mathbb{P}$. $\boldsymbol{\mu}$ and $\Sigma$ are the estimates of the mean and covariance matrix of the random vector $\boldsymbol{\lambda}$, respectively. $\mathcal{D}_{M1}(\mathcal{S}, \boldsymbol{\mu}, \Sigma, \gamma_1, \gamma_2)$ will also be referred to in short-hand notation as $\mathcal{D}_{M1}$ in this paper. Parameters $\gamma_1 \geq 0$ and $\gamma_2 \geq 0$ are responsible for controlling the "volume" of the ambiguity set.

*1) Data-Driven Uncertain Set $\mathcal{S}$*

Given a set $\{\boldsymbol{\lambda}^j\}_{j=1}^M$ of $M$ samples, now we are going to give the detailed data-driven ambiguity set.

Let $\lambda_{i,t}^- = \min\{\lambda_{i,t}^j, j = 1, \dots, M\}$, $\lambda_{i,t}^+ = \min\{\lambda_{i,t}^j, j = 1, \dots, M\}$, $\boldsymbol{\lambda}^- = [\lambda_{[i],[t]}^-]^T$, $\boldsymbol{\lambda}^+ = [\lambda_{[i],[t]}^+]^T$. In this paper, we use an uncertainty set to specify the lower bound $\lambda_{i,t}^-$ and upper bound $\lambda_{i,t}^+$ of each random variable $\lambda_{i,t}$, i.e.,

$$\hat{\mathcal{S}} = \{\boldsymbol{\lambda} \in \Re^{GT} | \boldsymbol{\lambda}^- \leq \boldsymbol{\lambda} \leq \boldsymbol{\lambda}^+\}. \tag{20}$$

(20) can be reformulated as:

$$\hat{\mathcal{S}} = \{\boldsymbol{\lambda} | \boldsymbol{A}\boldsymbol{\lambda} \leq \boldsymbol{B}\}, \quad (21)$$

where $\boldsymbol{A} = [\boldsymbol{I}; -\boldsymbol{I}]$, $\boldsymbol{B} = [\boldsymbol{\lambda}^+; -\boldsymbol{\lambda}^-]$.

*2) Data-Driven Ambiguity Set $\mathcal{D}_{M1}$*

We adopt a simple unbiased moment estimator [31]:

$$\hat{\boldsymbol{\mu}} = \frac{1}{M}\sum_{i=1}^{M}\boldsymbol{\lambda}_i, \quad (22)$$

$$\hat{\boldsymbol{\Sigma}} = \frac{1}{M}\sum_{i=1}^{M}(\boldsymbol{\lambda}_i - \hat{\boldsymbol{\mu}})(\boldsymbol{\lambda}_i - \hat{\boldsymbol{\mu}})^T. \quad (23)$$

Although $\hat{\boldsymbol{\Sigma}}$ is a positive semidefinite matrix by construction, we assume that $\hat{\boldsymbol{\Sigma}} \succ 0$ (in the case we study, we check this condition and it always holds unless $\boldsymbol{\lambda}_i = \hat{\boldsymbol{\mu}}, \forall j = 1, \ldots, M$).

Giving $\delta \in (0,1)$ and defining parameters $\bar{\delta}, \hat{R}, \bar{R}, \bar{a}, \bar{b}, \widehat{M}, \bar{\gamma}_1, \bar{\gamma}_2$ according to Eq. (24), we have that, if $M > \widehat{M}$ and $\boldsymbol{\lambda} \in \hat{\mathcal{S}}$, then with probability greater $1 - \delta$ the distribution of $\boldsymbol{\lambda}$ lies in the set $\mathcal{D}_{M1}(\hat{\mathcal{S}}, \hat{\boldsymbol{\mu}}, \hat{\boldsymbol{\Sigma}}, \bar{\gamma}_1, \bar{\gamma}_2)$. For more details, refer to Corollary 4 in [31].

$$\begin{cases} \bar{\delta} = 1 - \sqrt{1-\delta} \\ \hat{R} = \max_{\boldsymbol{\lambda} \in \hat{\mathcal{S}}} \|\hat{\boldsymbol{\Sigma}}^{-1/2}(\boldsymbol{\lambda} - \hat{\boldsymbol{\mu}})\| \\ \bar{R} = \hat{R}\left(1 - (\hat{R}^2 + 2)\left((2 + \sqrt{2\ln(4/\bar{\delta})})/\sqrt{M}\right)\right)^{-1/2} \\ \bar{a} = (\bar{R}^2/\sqrt{M})\left(\sqrt{1 - NT/\bar{R}^4} + \sqrt{\ln(4/\bar{\delta})}\right) \\ \bar{b} = (\bar{R}^2/\sqrt{M})\left(2 + \sqrt{2\ln(2/\bar{\delta})}\right)^2 \\ \widehat{M} = \max\left\{(\hat{R}^2 + 2)^2\left(2 + \sqrt{2\ln(4/\bar{\delta})}\right)^2, \frac{(8 + \sqrt{32\ln(4/\bar{\delta})})^2}{(\sqrt{\hat{R}+4}-\hat{R})^4}\right\} \\ \bar{\gamma}_1 = \frac{\bar{b}}{1-\bar{a}-\bar{b}} \\ \bar{\gamma}_2 = \frac{1+\bar{b}}{1-\bar{a}-\bar{b}} \end{cases} \quad (24)$$

Let $\boldsymbol{\zeta} = \hat{\boldsymbol{\Sigma}}^{-1/2}(\boldsymbol{\lambda} - \hat{\boldsymbol{\mu}})$, then (24)$^2$ (represents the 2$^{nd}$ formula of (24)) can be reformulated as

$$\hat{R} = \begin{cases} \max\|\boldsymbol{\zeta}\| \\ \text{s.t.} \hat{\boldsymbol{\Sigma}}^{-1/2}(\boldsymbol{\lambda}^- - \hat{\boldsymbol{\mu}}) \leq \boldsymbol{\zeta} \leq \hat{\boldsymbol{\Sigma}}^{-1/2}(\boldsymbol{\lambda}^+ - \hat{\boldsymbol{\mu}}) \end{cases}. \quad (25)$$

The close form solutions of (25) can be obtained as

$$\hat{R} = \max\{|\hat{\boldsymbol{\Sigma}}^{-1/2}(\boldsymbol{\lambda}^- - \hat{\boldsymbol{\mu}})|, |\hat{\boldsymbol{\Sigma}}^{-1/2}(\boldsymbol{\lambda}^+ - \hat{\boldsymbol{\mu}})|\}. \quad (26)$$

*D. Reformulation of DR-CVaR Model*

According to Lemma 1 of [31], by using duality theory in infinite dimensional convex problems and conic linear moment problems [32], [22], (18) is equivalent to the following convex semi-infinite optimization problem:

$$\min_{y,\alpha,Q,q,r,t} r + t \quad (27)$$

$$\text{s.t.} \begin{cases} r \geq h(\boldsymbol{y}, \alpha, \boldsymbol{\lambda}) - \boldsymbol{\lambda}^T \boldsymbol{Q}\boldsymbol{\lambda} - \boldsymbol{\lambda}^T \boldsymbol{q}, \forall \boldsymbol{\lambda} \in \hat{\mathcal{S}} \\ t \geq (\bar{\gamma}_2 \hat{\boldsymbol{\Sigma}} + \hat{\boldsymbol{\mu}}\hat{\boldsymbol{\mu}}^T) \circ \boldsymbol{Q} + \hat{\boldsymbol{\mu}}^T \boldsymbol{q} + \sqrt{\bar{\gamma}_1}\|\hat{\boldsymbol{\Sigma}}^{1/2}(\boldsymbol{q} + 2\boldsymbol{Q}\hat{\boldsymbol{\mu}})\|, \\ \boldsymbol{Q} \succeq 0, \boldsymbol{y} \in Y \end{cases}$$

*1) Linear Matrix Inequation (LMI) Reformulation*

(27)$^2$ is equivalent to

$$\begin{cases} r \geq \alpha - \boldsymbol{\lambda}^T \boldsymbol{Q}\boldsymbol{\lambda} - \boldsymbol{\lambda}^T \boldsymbol{q}, \forall \boldsymbol{\lambda} \in \hat{\mathcal{S}}, \\ r \geq \frac{1}{1-\beta}L(\boldsymbol{y},\boldsymbol{\lambda}) - \frac{\beta\alpha}{1-\beta} - \boldsymbol{\lambda}^T \boldsymbol{Q}\boldsymbol{\lambda} - \boldsymbol{\lambda}^T \boldsymbol{q}, \forall \boldsymbol{\lambda} \in \hat{\mathcal{S}}. \end{cases} \quad (28)$$

Let $g(\boldsymbol{\lambda}) = \boldsymbol{\lambda}^T \boldsymbol{Q}\boldsymbol{\lambda} + \boldsymbol{\lambda}^T \boldsymbol{q} + r - \alpha$, then (28)$^1$ is equivalent to

$$\min_{\boldsymbol{\lambda} \in \hat{\mathcal{S}}} g(\boldsymbol{\lambda}) \geq 0. \quad (29)$$

The Lagrangian dual problem of $\min_{\boldsymbol{\lambda} \in \hat{\mathcal{S}}} g(\boldsymbol{\lambda})$ is

$$\max_{\boldsymbol{\tau}_1^T \geq 0} \inf_{\boldsymbol{\lambda}} \{g(\boldsymbol{\lambda}) + \boldsymbol{\tau}_1^T(\boldsymbol{A}\boldsymbol{\lambda} - \boldsymbol{B})\}. \quad (30)$$

It is not difficult to verify that $\min_{\boldsymbol{\lambda} \in \hat{\mathcal{S}}} g(\boldsymbol{\lambda})$ is a convex programming problem and Slater's constraint qualification [34] holds for this problem. Then strong duality holds for $\min_{\boldsymbol{\lambda} \in \hat{\mathcal{S}}} g(\boldsymbol{\lambda})$ and (30). Furthermore, constraint (29) is equivalent to

$$\max_{\boldsymbol{\tau}_1^T \geq 0} \inf_{\boldsymbol{\lambda}} \{g(\boldsymbol{\lambda}) + \boldsymbol{\tau}_1^T(\boldsymbol{A}\boldsymbol{\lambda} - \boldsymbol{B})\} \geq 0, \quad (31)$$

(31) is further equivalent to the following constraint (the dimension of the vector is easily known from the context and will not be marked below):

$$\exists \boldsymbol{\tau}_1 \in \Re^{2GT} \geq 0, \boldsymbol{\lambda}^T \boldsymbol{Q}\boldsymbol{\lambda} + \boldsymbol{\lambda}^T \boldsymbol{q} + r - \alpha \geq -\boldsymbol{\tau}_1^T(\boldsymbol{A}\boldsymbol{\lambda} - \boldsymbol{B}), \forall \boldsymbol{\lambda}, \quad (32)$$

After applying the principles of Schur's complement, (32) is equivalent to (33).

$$\exists \boldsymbol{\tau}_1 \geq 0, \begin{bmatrix} \boldsymbol{Q} & \frac{1}{2}\boldsymbol{q} \\ \frac{1}{2}\boldsymbol{q}^T & r - \alpha \end{bmatrix} \succeq -\frac{1}{2}\begin{bmatrix} 0 & \boldsymbol{A}^T\boldsymbol{\tau}_1 \\ \boldsymbol{\tau}_1^T\boldsymbol{A} & -2\boldsymbol{\tau}_1^T\boldsymbol{B} \end{bmatrix}. \quad (33)$$

Similarly, (28)$^2$ is equivalent to

$$\exists \boldsymbol{\tau}_2 \geq 0, \begin{bmatrix} \boldsymbol{Q} & \frac{1}{2}\boldsymbol{q} \\ \frac{1}{2}\boldsymbol{q}^T & r \end{bmatrix} + \begin{bmatrix} 0 & \frac{1}{2}\frac{\boldsymbol{P}}{1-\beta} \\ \frac{1}{2}\frac{\boldsymbol{P}^T}{1-\beta} & \frac{\beta\alpha - \boldsymbol{1}^T\boldsymbol{Z}}{1-\beta} \end{bmatrix} \succeq -\frac{1}{2}\begin{bmatrix} 0 & \boldsymbol{A}^T\boldsymbol{\tau}_2 \\ \boldsymbol{\tau}_2^T\boldsymbol{A} & -2\boldsymbol{\tau}_2^T\boldsymbol{B} \end{bmatrix}. \quad (34)$$

*2) SDP Counterpart of DR-CVaR Self-Scheduling*

According to the LMI reformulations, (27) is equivalent to the SDP

$$\min_{\boldsymbol{y},\alpha,\boldsymbol{Q},\boldsymbol{q},r,t,\boldsymbol{\tau}_1 \geq 0, \boldsymbol{\tau}_2 \geq 0} r + t, \quad (35)$$

$$\text{s.t.} \begin{cases} (33)(34) \\ t \geq (\bar{\gamma}_2\hat{\boldsymbol{\Sigma}} + \hat{\boldsymbol{\mu}}\hat{\boldsymbol{\mu}}^T) \circ \boldsymbol{Q} + \hat{\boldsymbol{\mu}}^T\boldsymbol{q} + \sqrt{\bar{\gamma}_1}\|\hat{\boldsymbol{\Sigma}}^{1/2}(\boldsymbol{q} + 2\boldsymbol{Q}\hat{\boldsymbol{\mu}})\| \\ \boldsymbol{Q} \succeq 0, \boldsymbol{y} \in Y \end{cases}$$

## IV. APPROXIMATIONS OF DR-CVaR MODEL

Although SDP formulations (35) is polynomially solvable in theory, it require significant computational efforts because of the high-dimensional matrix constraints of $\boldsymbol{Q}$. To overcome such challenges, in this section, we will introduce two approximations with smaller-size SDP matrix constraints.

*A. Upper Approximation Based on Vector Splitting*

First, we perform an eigenvalue decomposition on matrix $\boldsymbol{\Sigma} = U\Lambda U^T = U\Lambda^{1/2}(U\Lambda^{1/2})^T$. Then, let

$$\boldsymbol{\lambda}_c = (U\Lambda^{-1/2})^T(\boldsymbol{\lambda} - \boldsymbol{\mu}), \quad (36)$$

We have

$$\boldsymbol{\lambda} = U\Lambda^{1/2}\boldsymbol{\lambda}_c + \boldsymbol{\mu}, \quad (37)$$

$$\mathcal{S}_c(\mathcal{S}) = \{\boldsymbol{\lambda}_c | U\Lambda^{1/2}\boldsymbol{\lambda}_c + \boldsymbol{\mu} \in \mathcal{S}\}. \quad (38)$$

We derive the approximation by splitting the random vector $\boldsymbol{\lambda}_c$ into $P$ pieces, i.e., $\boldsymbol{\lambda}_c = (\boldsymbol{\lambda}_{c_1}; \boldsymbol{\lambda}_{c_2}; \ldots; \boldsymbol{\lambda}_{c_P})$, where $\boldsymbol{\lambda}_{c_i} \in \Re^{m_i}, \forall i \in [P]$ and $\sum_{i=1}^{P} m_i = GT$. After ignoring covariances (some correlations) among $\boldsymbol{\lambda}_{c_p}$ and $\boldsymbol{\lambda}_{c_q}$ for any $p, q \in [P]$ with $p \neq q$. We obtain the following ambiguity set:

$$\mathcal{D}_{M2}(\mathcal{S}_c(\mathcal{S}), \boldsymbol{\mu}, \boldsymbol{\Sigma}, \gamma_1, \gamma_2, P) = \left\{ \mathbb{P}_c \in \mathcal{P}_0(\Re^{GT}) \middle| \begin{matrix} \mathbb{P}(\boldsymbol{\lambda}_c \in \mathcal{S}_c(\mathcal{S})) = 1 \\ \mathbb{E}_{\mathbb{P}_c}(\boldsymbol{\lambda}_c)^T \mathbb{E}_{\mathbb{P}_c}(\boldsymbol{\lambda}_c) \leq \gamma_1 \\ \mathbb{E}_{\mathbb{P}_c}[\boldsymbol{\lambda}_{c_i}\boldsymbol{\lambda}_{c_i}^T] \preceq \gamma_2 \boldsymbol{I}_{m_i}, \forall i \in [P] \end{matrix} \right\}. \quad (39)$$

It is not hard to verify that $\mathcal{D}_{M1}(\mathcal{S}, \boldsymbol{\mu}, \boldsymbol{\Sigma}, \gamma_1, \gamma_2)$ is equivalent to $\mathcal{D}_{M2}(\mathcal{S}_c(\mathcal{S}), \boldsymbol{\mu}, \boldsymbol{\Sigma}, \gamma_1, \gamma_2, 1)$, which is a subset of $\mathcal{D}_{M2}(\mathcal{S}_c(\mathcal{S}), \boldsymbol{\mu}, \boldsymbol{\Sigma}, \gamma_1, \gamma_2, P)$ [35]. Then according to the last subsection, if $M > \widehat{M}$ and $\boldsymbol{\lambda}_c \in \mathcal{S}_c(\hat{\mathcal{S}})$, with probability greater $1 - \delta$ the distribution of $\boldsymbol{\lambda}_c$ lies in the set $\mathcal{D}_{M2}(\mathcal{S}_c(\hat{\mathcal{S}}), \hat{\boldsymbol{\mu}}, \hat{\boldsymbol{\Sigma}}, \bar{\gamma}_1, \bar{\gamma}_2, P)$. Then, we obtain the upper ap-

proximation of (18):
$$\min_{\mathbf{y}\in Y,\alpha\in\Re}\left\{\sup_{\mathbb{P}_c\in\mathcal{D}_{M2}}\{\mathbb{E}_{\mathbb{P}}(h(\mathbf{y},\alpha,U\Lambda^{1/2}\boldsymbol{\lambda}_c+\widehat{\boldsymbol{\mu}}))\}\right\}, \quad (40)$$

Similar to the reformulation (18) of in last subsection, (40) is equivalent to the following convex optimization problem:
$$\min_{\mathbf{y},\alpha,Q_i,q,r,t} r + t, \quad (41)$$
$$\text{s.t.}\begin{cases} r \geq h(\mathbf{y},\alpha,U\Lambda^{1/2}\boldsymbol{\lambda}_c+\widehat{\boldsymbol{\mu}}) - \sum_{i=1}^P \boldsymbol{\lambda}_{c_i}^T Q_i \boldsymbol{\lambda}_{c_i} - q^T \boldsymbol{\lambda}_c, \forall \boldsymbol{\lambda}_c \in \mathcal{S}_c(\hat{S}) \\ t \geq \sum_{i=1}^P \{(\bar{\gamma}_2 \boldsymbol{I}_{m_i}) \circ \boldsymbol{Q}_i\} + \sqrt{\bar{\gamma}_1}\|q\| \\ \boldsymbol{Q}_i \succcurlyeq 0, \forall i \in [P] \\ \mathbf{y} \in Y \end{cases}.$$

$(41)^2$ is equivalent to
$$\begin{cases} r \geq \alpha - \sum_{i=1}^P \boldsymbol{\lambda}_{c_i}^T \boldsymbol{Q}_i \boldsymbol{\lambda}_{c_i} - \boldsymbol{\lambda}_c^T q, \forall \boldsymbol{\lambda}_c \in \mathcal{S}_c(\hat{S}), \\ r \geq \frac{1}{1-\beta} L(\mathbf{y},U\Lambda^{1/2}\boldsymbol{\lambda}_c+\widehat{\boldsymbol{\mu}}) - \frac{\beta\alpha}{1-\beta} - \sum_{i=1}^P \boldsymbol{\lambda}_{c_i}^T \boldsymbol{Q}_i \boldsymbol{\lambda}_{c_i} - \boldsymbol{\lambda}_c^T q, \forall \boldsymbol{\lambda}_c \in \mathcal{S}_c(\hat{S}). \end{cases} \quad (42)$$

Similar to the variation of $(28)^1$, $(42)^1$ is equivalent to:
$$\exists \tau_1 \geq 0, \sum_{c_i=1}^p \boldsymbol{\lambda}_{c_i}^T \boldsymbol{Q}_i \boldsymbol{\lambda}_{c_i} + \boldsymbol{\lambda}_c^T q + r - \alpha + \tau_1^T (\boldsymbol{A}(U\Lambda^{1/2}\boldsymbol{\lambda}_c + \widehat{\boldsymbol{\mu}}) - \boldsymbol{B}) \geq 0, \forall \boldsymbol{\lambda}_{c_i}, \quad (43)$$

Then, we perform the following decomposition:
$$U\Lambda^{1/2}\boldsymbol{\lambda}_c = \sum_{i=1}^p U_{GT,m_i}\Lambda_{m_i}^{1/2}\boldsymbol{\lambda}_{c_i}, \quad (44)$$
where $U_{GT,m_i}$ is $[1:GT, \sum_{j=1}^i m_{i-1}+1:\sum_{j=1}^i m_i]$ submatrix. $\Lambda_{m_i}^{1/2}$ is upper-left submatrix of $\Lambda^{1/2}$.

By plugging (44) to (43), we have
$$\exists \tau_1 \geq 0, r - \alpha - \tau_1^T \boldsymbol{B} + \tau_1^T A\widehat{\boldsymbol{\mu}} + \sum_{i=1}^p \left(\boldsymbol{\lambda}_{c_i}^T \boldsymbol{Q}_i \boldsymbol{\lambda}_{c_i} + \left(q_{c_i}^T + \tau_1^T \boldsymbol{A} U_{GT,m_i}\Lambda_{m_i}^{1/2}\right)\boldsymbol{\lambda}_{c_i}\right) \geq 0, \forall \boldsymbol{\lambda}_{c_i}, \quad (45)$$

By introducing auxiliary variables, we have
$$\exists \tau_1 \geq 0, \begin{cases} r - \alpha - \tau_1^T \boldsymbol{B} + \tau_1^T A\widehat{\boldsymbol{\mu}} + \sum_{i=1}^n z_i^1 \geq 0 \\ \boldsymbol{\lambda}_{c_i}^T \boldsymbol{Q}_i \boldsymbol{\lambda}_{c_i} + (q_{m_i}^T + \tau_1^T \boldsymbol{A} U_{GT,m_i}\Lambda_{m_i}^{1/2})\boldsymbol{\lambda}_{c_i} \geq z_i^1 \end{cases}, \forall \boldsymbol{\lambda}_{c_i}, \quad (46)$$

After applying the principles of Schur's complement, (46) is equivalent to
$$\exists \tau_1 \geq 0, \begin{cases} r - \alpha - \tau_1^T \boldsymbol{B} + \tau_1^T A\widehat{\boldsymbol{\mu}} + \sum_{i=1}^P z_i^1 = 0 \\ \begin{bmatrix} \boldsymbol{Q}_i & \frac{1}{2}\omega_1(\tau_1) \\ \frac{1}{2}(\omega_1(\tau_1))^T & -z_i^1 \end{bmatrix} \succcurlyeq 0, \forall i \in [P] \end{cases}, \quad (47)$$
where $\omega_1(\tau_1) = (q_{m_i}^T + \tau_1^T \boldsymbol{A} U_{GT,m_i}\Lambda_{m_i}^{1/2})^T$.

Similarly, $(42)^2$ is equivalent to
$$\exists \tau_2 \geq 0$$
$$\begin{cases} r - \tau_2^T \boldsymbol{B} + \tau_2^T A\widehat{\boldsymbol{\mu}} + \frac{\beta\alpha}{1-\beta} - \frac{1}{1-\beta}\boldsymbol{1}^T z + \frac{1}{1-\beta}\boldsymbol{P}^T\widehat{\boldsymbol{\mu}} + \sum_{i=1}^P z_i^2 = 0 \\ \begin{bmatrix} \boldsymbol{Q}_i & \frac{1}{2}\omega_2(\tau_2) \\ \frac{1}{2}(\omega_2(\tau_2))^T & -z_i^2 \end{bmatrix} \succcurlyeq 0, \forall i \in [P] \end{cases}, \quad (48)$$
where $\omega_2(\tau_2) = \omega_1(\tau_2) + (\frac{\boldsymbol{P}^T U_{GT,m_i}\Lambda_{m_i}^{1/2}}{1-\beta})^T$.

Then, (41) is equivalent to the SDP with smaller-size matrix constraints
$$\min_{\mathbf{y},\alpha,Q_i,q,r,t,\tau_1\geq 0,\tau_2\geq 0,z_i^1,z_i^2} r + t, \quad (49)$$
$$\text{s.t.}\begin{cases} (47)(48) \\ t \geq \sum_{i=1}^P \{(\bar{\gamma}_2 \boldsymbol{I}_{m_i}) \circ \boldsymbol{Q}_i\} + \sqrt{\bar{\gamma}_1}\|q\| \\ \boldsymbol{Q}_i \succcurlyeq 0, \forall i \in [P] \\ \mathbf{y} \in Y \end{cases}.$$

### B. Approximation Based on Region Partition

Similar to the last subsection, we derive another approximation by splitting the random vector $\boldsymbol{\lambda}$ into $P$ pieces, i.e., $\boldsymbol{\lambda} = (\boldsymbol{\lambda}_1; \boldsymbol{\lambda}_2; \dots; \boldsymbol{\lambda}_P)$, where $\boldsymbol{\lambda}_d \in \Re^{m_d}, \forall d \in [P]$ and $\sum_{d=1}^P m_d = GT$. $\mathbf{y}_d \in Y_d, \mathbf{y} \in \bar{Y}, (\bigcap_{d\in[P]} Y_d) \cap \bar{Y} = Y$. For the $d$th region, we add lower-script "$d$" to the variables and functions which corresponding to this region.

We assume that $\boldsymbol{\lambda}_d \sim \mathbb{P}_d$ (which is an approximation for $\boldsymbol{\lambda} \sim \mathbb{P}$), i.e., we ignore all correlations of $\boldsymbol{\lambda}_d$ among different regions. Then, similar to $\mathcal{D}_{M1}(\hat{S}, \widehat{\boldsymbol{\mu}}, \hat{\Sigma}, \bar{\gamma}_1, \bar{\gamma}_2)$, we can separately construct the ambiguity set for $\mathbb{P}_d$, which can be denoted as $\mathcal{D}_{M1,d}(\hat{S}_d, \widehat{\boldsymbol{\mu}}_d, \hat{\Sigma}_d, \bar{\gamma}_{1,d}, \bar{\gamma}_{2,d})$. We obtain the following approximation of $\mathcal{D}_{M1}$:
$$\mathcal{D}_{M3} = \prod_{d=1}^P \mathcal{D}_{M1,d}. \quad (50)$$

Then our approximation of DR-CVaR model of self-scheduling based on region partition is
$$\min_{\mathbf{y}_d\in Y_d,\mathbf{y}\in\bar{Y},\alpha_d\in\Re} \sum_{d\in[P]}\left\{\sup_{\mathbb{P}_d\in\mathcal{D}_{M1,d}}\{\mathbb{E}_{\mathbb{P}}(h_d(\mathbf{y}_d,\alpha_d,\boldsymbol{\lambda}_d))\}\right\}. \quad (51)$$

Same as above, (51) is equivalent to the SDP
$$\min_{\mathbf{y}_d,\alpha_d,Q_d,q_d,r_d,t_d,\tau_{1,d}\geq 0,\tau_{2,d}\geq 0} \sum_{d\in[P]}\{r_d + t_d\}, \quad (52)$$
$$\text{s.t.}\begin{cases} \begin{bmatrix} \boldsymbol{Q}_d & \frac{1}{2}q_d \\ \frac{1}{2}q_d^T & r_d - \alpha_d \end{bmatrix} \succcurlyeq -\frac{1}{2}\begin{bmatrix} 0 & \boldsymbol{A}_d^T\tau_{1,d} \\ \tau_{1,d}^T \boldsymbol{A}_d & -2\tau_{1,d}^T \boldsymbol{B}_d \end{bmatrix} \\ \begin{bmatrix} \boldsymbol{Q}_d & \frac{1}{2}q_d \\ \frac{1}{2}q_d^T & r_d - \alpha_d \end{bmatrix} + \begin{bmatrix} 0 & -\frac{1}{2}\frac{\boldsymbol{P}_d}{1-\beta_d} \\ -\frac{1}{2}\frac{\boldsymbol{P}_d^T}{1-\beta_d} & \frac{\alpha_d}{1-\beta_d} + f_d(\boldsymbol{P}_d) \end{bmatrix} \succcurlyeq -\frac{1}{2}\begin{bmatrix} 0 & \boldsymbol{A}_d^T\tau_{2,d} \\ \tau_{2,d}^T \boldsymbol{A}_d & -2\tau_{2,d}^T \boldsymbol{B}_d \end{bmatrix} \\ t_d \geq (\bar{\gamma}_{2,d}\hat{\Sigma}_d + \widehat{\boldsymbol{\mu}}_d\widehat{\boldsymbol{\mu}}_d^T) \circ \boldsymbol{Q}_d + \widehat{\boldsymbol{\mu}}_d q_d + \sqrt{\bar{\gamma}_{1,d}}\|\hat{\Sigma}_d^{1/2}(q_d + 2\boldsymbol{Q}_d\widehat{\boldsymbol{\mu}}_d)\| \\ \boldsymbol{Q}_d \succcurlyeq 0 \\ \mathbf{y}_d \in Y_d, \mathbf{y} \in \bar{Y} \end{cases}.$$

We note that the SDP (52), similar to (49), only includes smaller-size matrix constraints. Furthermore, except for $\mathbf{y} \in \bar{Y}$, all the other parts of this problem are separable in the regions. The ADMM can be used to solve this problem distributionally.

## V. SIMULATION RESULT

For verify the correctness and economy of the proposed model, IEEE 6-bus, 30-bus and 118-bus test system are selected for simulation calculation. The calculation environment is Win10 system, CPU is Inter Core i5-4210M with 2.60GHz CPU clock speed. The model is solved by calling CPLEX and MOSEK from Yalmip. The generator parameters of 6-bus and 30-bus are given in Table1 and Table 2 respectively. Due to excessive amount data of 118-bus system, for the sake of brevity, it is not listed here [36].

TABLE 1
Generator data of 6-bus test system

| Bus No. | $\bar{P}_i^G$[MW] | $\underline{P}_i^G$[MW] | Cost coefficients | | |
|---|---|---|---|---|---|
| | | | a[$/h] | b[$/MW h] | c[$/MWh$^2$] |
| 1 | 200 | 50 | 0.00 | 2.00 | 0.00375 |
| 2 | 80 | 20 | 0.00 | 1.75 | 0.01750 |
| 6 | 50 | 15 | 0.00 | 1.00 | 0.06250 |

TABLE 2
Generator data of 30-bus test system

| Bus No. | $\bar{P}_i^G$[MW] | $\underline{P}_i^G$[MW] | Cost coefficients | | |
|---|---|---|---|---|---|
| | | | a[$/h] | b[$/MW h] | c[$/MWh$^2$] |
| 1 | 60 | 10 | 10 | 0.020 | 0.00020 |
| 2 | 60 | 10 | 10 | 0.015 | 0.00024 |
| 5 | 150 | 10 | 20 | 0.018 | 0.00008 |
| 8 | 120 | 10 | 10 | 0.010 | 0.00012 |
| 11 | 150 | 10 | 20 | 0.018 | 0.00008 |
| 13 | 60 | 10 | 10 | 0.015 | 0.00020 |

### A. Simulation of uncertain electricity price

Firstly, the quadratic programming problem with the objective function of generating cost considering unit operating constraints and network power flow constraints is solved to obtain the reasonable price. Then, in order to be closer to the uncertainty of the actual scenario, the electricity price data set used in following experiment is obtained by random fluctuation of the price. The electricity price fluctuations of unit1, unit2 and unit3 of the IEEE 6-bus 24-period test system are shown in Fig.1. (There are 168 samples in each period in this case.)

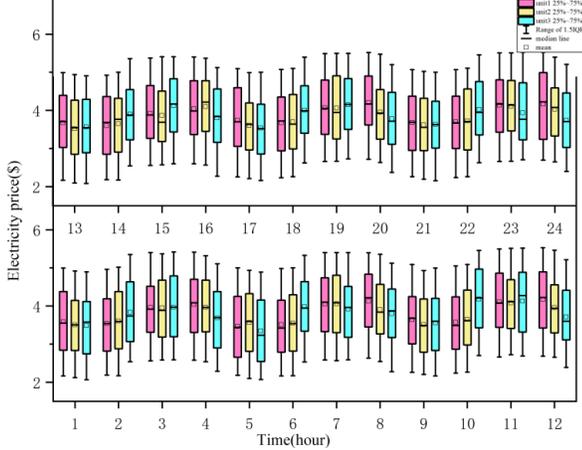

Fig.1 Box-plot of electricity price

### B. Verification of model characteristics

In this subsection, IEEE 6-bus 4-periods test system is selected for stimulations. Table 3 shows the change of profit and object value with the change of ambiguity set, where the profit is calculated by formulate $F(y,\lambda) = P^T\lambda - \mathbf{1}^T z$, and the object value denotes $\beta$-CVaR, that is, risk. All the codes and cases of the simulations for this paper can be freely downloaded from https://github.com/linfengYang/A-novel-DRO-model-for-self-scheduling-problem.

TABLE 3
The change of profit and risk in different ambiguity set

| | | | | |
|---|---|---|---|---|
| | | $\delta=0.10$ | | |
| $\beta$ | 0.85 | 0.90 | 0.95 | 0.99 |
| object | 290.70 | 209.10 | 76.74 | 5.51 |
| profit | 702.29 | 601.14 | 506.00 | 499.74 |
| | | $\delta=0.20$ | | |
| $\beta$ | 0.85 | 0.90 | 0.95 | 0.99 |
| object | 294.87 | 212.55 | 79.13 | 5.51 |
| profit | 719.21 | 601.95 | 507.19 | 499.74 |
| | | $\delta=0.30$ | | |
| $\beta$ | 0.85 | 0.90 | 0.95 | 0.99 |
| object | 297.64 | 214.78 | 80.65 | 5.51 |
| profit | 729.14 | 602.50 | 508.50 | 499.74 |
| | | $\delta=0.40$ | | |
| $\beta$ | 0.85 | 0.90 | 0.95 | 0.99 |
| object | 299.84 | 216.51 | 81.84 | 5.51 |
| profit | 737.13 | 602.92 | 509.14 | 499.74 |

DRO can weigh the robustness and economy by controlling the size of ambiguity set. It can be seen from formula (24) that the ambiguity set increases with the decrease of $\delta$, and combined with Table 3, we can find that the larger the ambiguity set is, the more probability distribution it contains, and the more conservative strategy is obtained. In this case, for the sake of the possible risks, each unit generates less power, so the profit is correspondingly reduced. Users can choose different size of ambiguity sets according to different scenario, so as to choose whether they prefer robustness or economy.

### C. Contrast with the RO model

To consider both robustness and economy, we fixed $\delta$ to 0.2 in following experiments. At the same time, considering the practical significance of $\beta$, set $\beta$ as 0.90, 0.95 and 0.99. The RO model is implemented refer to [1], and IEEE 6- and 30-bus power system was used for the control experiment. Fig.2 shows the relationship between total output and load of the DRO model and RO model under 6-24 test system with $\beta$ fixed at 0.90, the output of unit1, unit2 and unit3 in one day (24 hours) in different models of this situation is shown in Fig.3. And Fig.4 shows the line chart of two models' profit changing with the $\beta$ under 6-24 and 30-24 test system. Finally, the comparison of the calculation time of the two models in different test system and different sample numbers is shown in Fig.5.

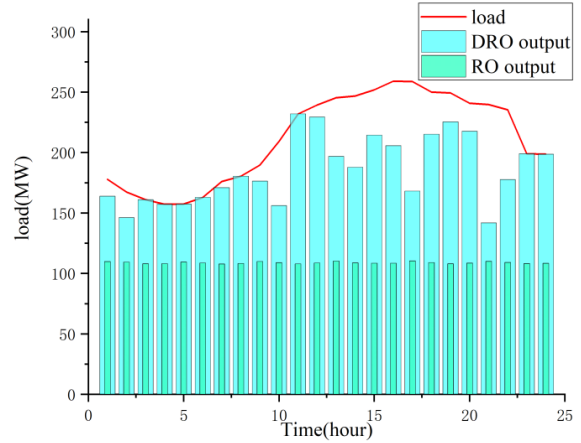

Fig.2 Diagram of total output verses load of two models

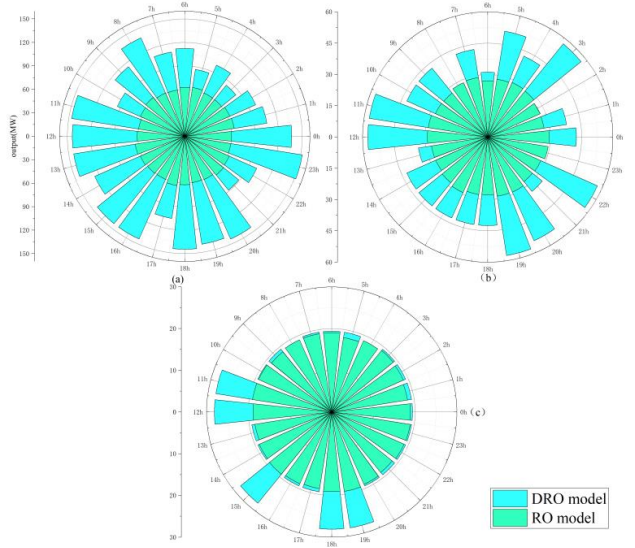

Fig.3 Units output comparison diagram of two model: (a) unit 1; (b) unit 2; (c) unit 3

It can be seen from Fig.2 and Fig.3 that compared with the DRO model, the RO model overemphasizes the robustness for make the potential loss as small as possible when the electricity price fluctuates, which results in that the units in the RO model are unwilling to generate electricity, so the total generating capacity in the obtained RO strategy is usually only about 2/3 of the load requirement, and the scheduling strategy basically does not change with the change of load requirement, but fluctuates slightly near the total generating capacity which can most resist the risk. And combined with Fig.1, it is not difficult to find that the generation of RO model units is less sensitive to fluctuations of electricity price in various periods.

However, each unit of the DRO model is more sensitive to the fluctuations of electricity price in each period, and the power generation of the unit will change with the price of electricity, moreover, combined with the cost coefficients of the units in Table 1, the power generation of each unit not only considers the fluctuations in electricity price but also relates to its own cost coefficient (For example, at 11:00, the electricity price of all three units is high and the variance is small, so all three units generate a lot of electricity at this time; at 16:00, although the load requirement is higher, unit 2 and unit 3 do not vigorously generate electricity due to large variance of price and low price respectively.) Among them, the electricity price fluctuation of unit 1 is small and the quadratic term of the cost coefficient also small, so unit 1 undertakes the main power generation task, and its power generation curve basically conforms to the load curve, and the electricity price fluctuations and quadratic term of the cost coefficient of unit 3 are largest, so considering the robustness, the power generation of unit 3 is stable except for periods when the electricity price is high, and unit 2 is in between.

Through the above numerical experimental results, we can find that compared to the RO model, the DRO model we proposed more comprehensively considers electricity price fluctuations, unit parameters and load requirement. Therefore, it can increase the economy while ensuring the required robustness. These characteristics will be proved in the following experimental analysis.

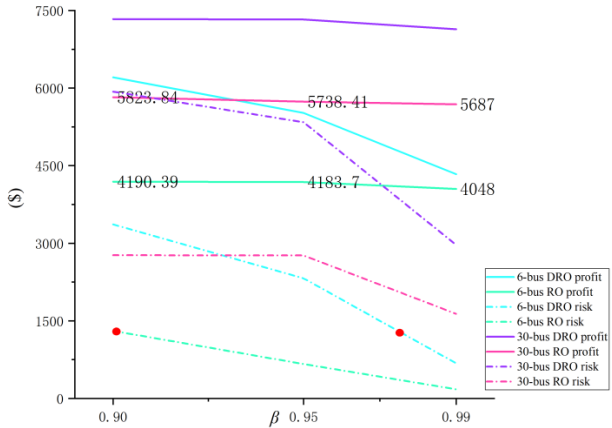

Fig.4 Profit and risk comparison chart of two models

Combined with Fig.4, it can be concluded that RO model can indeed avoid potential risk better than DRO model, however, due to its too conservative self-scheduling strategy, it cannot mobilize the enthusiasm of generating units. Furthermore, because it cares too much about robustness, it ignores economy, the expected profit of the RO scheduling strategy is usually less than that of DRO. At the same time, DRO model does not blindly ignore the robustness and pursue the economy (As shown by the two red dots in Fig.4, the risks of the two models are the same and the profit of DRO model is clearly higher at the red dots, therefore, DRO model can improve profit on the premise of ensuring risk management.) We can see from the analysis in the previous paragraph that although the potential risks of DRO model are relatively large and some robustness is given up for economic considerations, its scheduling strategy is to comprehensively consider the cost coefficient of each unit, the fluctuation of electricity prices in each period, and load requirements, etc. So it is concluded that it is completely reasonable to reduce the weight of robustness on this basis.

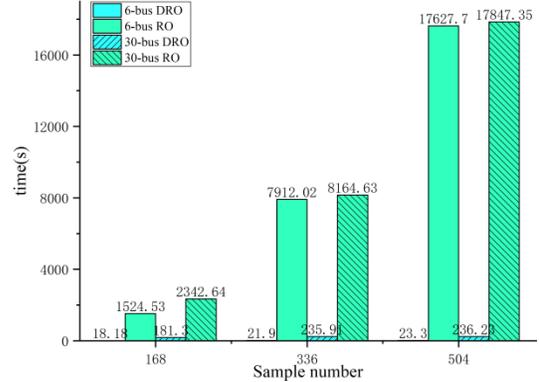

Fig.5 Calculation time comparison chart of two models

Then through Fig.5 and the formula of RO model [1], we can find that with the increase of the number of uncertainty variable samples, time periods and units, the solution time of RO model will increase sharply. However, the more samples the RO model has, the accurate it can be, this is contrary to the above characteristics. At the same time, the RO model is difficult to deal with large-scale systems. Nevertheless, when deal with the uncertainty variable samples, the DRO model calculates the moment information about the uncertainty variable, and then uses the moment information to construct the ambiguity set, so that when the number of samples and periods increases, the DRO model can still deal with the uncertainty variable quickly.

### D. Simulation result of approximate models

Firstly, the DRO model and two approximate models (one approximate model based on vector splitting called APP1, and another based on ADMM called APP 2.) is applied on IEEE 30-bus test system to illustrate the approximate models' accuracy. In fact, APP2 can also be solved by directly calling solver, then the computational time and resources required for the two approximate models are roughly same, but in this case, the property of block separation is not well utilized. Therefore, the ADMM algorithm is implemented to solve APP2 by referring to [37]. In this subsection, we refer to [38] to achieve the division of area, and fixed $\delta$ to 0.2 and $\beta$ to 0.90 to calculate the gap between the approximate models and the DRO model profit in the same situation, the gap is calculated by formula $gap = (|P_{DRO} - P_{APP}|)/P_{DRO}$, where $P(\cdot)$ is the profit. To ensure that the approximate models can approximate the DRO model stably and accurately, we generated 50 copies of electricity price data for simulation, and the simulation results will be displayed as a group of every five tests. Fig.6 (a) and (b) show the results of APP1 and APP2, respectively. Then, for verify the characteristics of APP2, 118-bus system is divided into 3 areas and 10 areas, the results of the solution time of the DRO model and approximate models on the IEEE 30-bus and 118-bus test system are shown in TABLE 4, to verify that the approximate models greatly shorten the calculation time and requires less computing resource. The DRO model cannot be solved effec-tively in this experimental environment because of the large scale of the SDP problem under 118-bus test system.

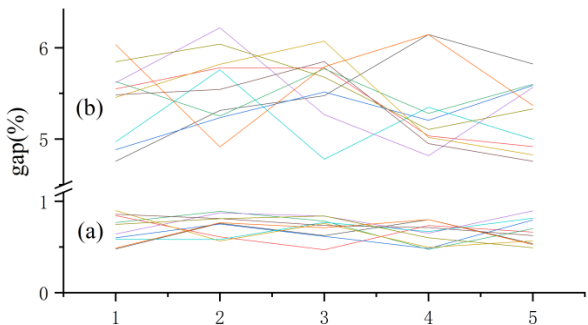

Fig.6. Line chart of gap: (a) Gap of the DRO model and APP1 model; (b) Gap of the DRO model and APP2 model

TABLE 4

|  |  | DRO | APP1 | APP2 | APP2 |
|---|---|---|---|---|---|
| 30-bus | areas | 1 | 2 | 2 | |
| | time | 231.64s | 35.88s | 123.35s | |
| 118-bus | areas | 1 | 6 | 3 | 10 |
| | time | | 3h16m | 10h51m | 2h18m |

As shown in Fig.6, the two approximate models can effectively approximate the DRO model, combined with Table 4, we can draw a further conclusion that APP1 can approximate the DRO model more accurately, its gap is always less than 0.9%, and the gap between DRO model and APP2 will be slightly larger. However, since APP2 can be calculated in parallel, the solution time will be greatly reduced when the number of areas increases, and APP2 uses a distributed method to solve, so it has greater value in practical applications, the units or power areas can be managed separately, they are no need to share their key information, e.g., types, characteristics, capacities of energy source and loads, to others for the sake of privacy.

## VI. CONCLUSION

This paper proposes a novel distributionally robust optimization model for generation companies. For self-scheduling problem in uncertain environment, an adjustable strategy is obtained that considers unit cost parameters, electricity price fluctuation and load requirement, which can handle effectively uncertainties associate with electricity price, even when its probability distribution are not known perfectly. Furthermore, two high quality approximate models are developed because of the difficulty of solving the large-scale SDP, which greatly reduce the time and resources required for calculation while ensuring the accuracy. Simulations are conducted to show the characteristics of proposed DRO model and the accuracy and efficiency of the two approximate models.